\documentclass[12pt]{article}
\usepackage{amssymb}

\newtheorem{theorem}{Theorem}

\newtheorem{lemma}[theorem]{Lemma}

\newtheorem{proposition}[theorem]{Proposition}

\begin{document}

\title{Ornstein-Uhlenbeck semi-groups on stratified groups}
\author{Fran\c{c}oise Lust-Piquard}
\maketitle

\begin{abstract}
We consider, in the setting of stratified groups $G,$ two analogues of the
Ornstein-Uhlenbeck semi-group, namely Markovian diffusion semi-groups acting
on $L^{q}(p(\gamma )d\gamma ),$ whose invariant density $p$ is a heat kernel
at time 1 on $G.$

The first one is symmetric on $L^{2}(pd\gamma ),$ its generator is $%
\sum_{i=1}^{n}X_{i}^{\ast }X_{i},$ where $(X_{i})_{i=1}^{n}$ \ is a basis of
the first layer of the Lie algebra of $G.$

The second one, denoted by $T_{t}=e^{-tN},t>0,$ is non symmetric on $%
L^{2}(pd\gamma )$ and the formal real part of $N$ is $\sum_{i=1}^{n}X_{i}^{%
\ast }X_{i}.$ The operators $e^{-tN}$ are compact on $L^{q}(pd\gamma ),$ $%
1<q<\infty .$ The spectrum of $N$ on this space is the set of integers $%
\mathbb{N}$ if polynomials are dense in $L^{2}(p(\gamma )d\gamma ),$ i.e if $%
G$ has at most 4 layers; and we determine in this case its eigenspaces. When 
$G$ is step 2, we give another description of these eigenspaces, very
similar to the classical definition of "Hermite polynomials" by their
generating function.


\end{abstract}

\section{\textbf{Introduction and notation}}

Let $G$ be a stratified Lie group equipped with its (biinvariant) Haar
measure $dg$ and dilations $(\delta _{t})_{t\geq 0}$. Let $Q$ be the
homogeneous dimension of $G$. We denote by $\mathcal{D}(G)$ the space of $%
\mathcal{C}^{\infty }$ compactly supported functions on $G,$ by $\mathcal{S}%
(G)$ the space of Schwartz functions, by $\mathcal{S}^{\prime }(G)$ its
dual, and $L^{q}(\varphi dg)=L^{q}(G,\varphi dg)$ for a measurable non
negative function $\varphi $.

\noindent As usual, elements $Z$ of the Lie algebra $\mathcal{G}$ are
identified with left invariant vector fields by

\[
(Zf)(g)=\frac{d}{dt}\mid _{t=0}f(g\exp tZ). 
\]

\noindent Let $L$ be a subLaplacian on $G,$ i.e. an operator on $\mathcal{S(}%
G\mathcal{)}$ defined by

\begin{equation}
L=-\sum_{1}^{n}X_{i}^{2}  \label{0}
\end{equation}%
where $(X_{i})_{1\leq i\leq n}$ is a linear basis of the first layer of $%
\mathcal{G}.$ Obviously $L$ commutes with left translations and satisfies

\begin{equation}
\delta _{t^{-1}}L\delta _{t}=t^{2}L,\;t>0.  \label{1}
\end{equation}

\noindent The following facts are well known, see e.g. \cite[propositions
1.68, 1.70, 1.74]{FS}: $-\frac{L}{2}$ generates a strongly continuous
semi-group $e^{-\frac{t}{2}L}$ of convolution operators which are
contractions on $L^{q}(dg),$ $1\leq q\leq \infty .$ The kernel $p_{t}$ of $%
e^{-\frac{t}{2}L}$ is a positive function such that $p_{t}(g)=p_{t}(g^{-1}),$
it\ lies in $\mathcal{S(}G\mathcal{)}$ and has norm one in $L^{1}(dg).$
Denoting $p_{1}=p,$

\[
p_{t}(g)=t^{-\frac{Q}{2}}p\circ \delta _{\frac{1}{\sqrt{t}}}(g). 
\]

\noindent Equivalently, for $f\in L^{q}(dg),$

\begin{equation}
e^{-\frac{t}{2}L}(f)(\gamma )=f\ast p_{t}(\gamma )=\int_{G}f(\gamma
g^{-1})p_{t}(g)dg=\int_{G}f(\gamma \delta _{\sqrt{t}}g^{-1})p(g)dg.
\label{2}
\end{equation}

The aim of this paper is to generalize the Ornstein-Uhlenbeck semi-group in
the setting of stratified groups, namely to consider Markovian semi-groups
acting on $L^{q}(p(\gamma )d\gamma ),1\leq q\leq \infty ,$ for which $%
p(\gamma )d\gamma $ is an invariant measure, whose generators are related to
the first layer gradient

\[
\nabla =(X_{1},..,X_{n}). 
\]

\noindent The classical Ornstein-Uhlenbeck semi-group is defined on $%
\mathcal{S}(\mathbb{R}^{n})$ by Mehler formula

\[
e^{-tN_{0}}(f)(x)=\int_{\mathbb{R}^{n}}f(e^{-t}x+\sqrt{1-e^{-2t}}%
y)p(y)dy,\;t\geq 0, 
\]

\noindent where the gaussian density $p(y)=\frac{1}{(2\pi )^{\frac{n}{2}}}%
e^{-\frac{1}{2}\left\vert y\right\vert ^{2}}$ is the kernel of $e^{-\frac{%
\Delta }{2}},$ and $\Delta $ is the (positive) Laplacian on $\mathbb{R}^{n}$%
. The O-U semi-group is contracting on $L^{q}(\mathbb{R}^{n},pdx),$ $1\leq
q\leq \infty ,$ compact if $1<q<\infty ,$ but not compact on $L^{1}(\mathbb{R%
},pdx)$ \cite[theorem 4.3.5]{D}, and $p$ is an invariant measure. The
generator $-N_{0}$ satisfies

\[
N_{0}=\sum_{j=1}^{n}(\frac{\partial }{\partial x_{j}})^{\ast }\frac{\partial 
}{\partial x_{j}}=\Delta -\sum_{j=1}^{n}\frac{\frac{\partial p}{\partial
x_{j}}}{p}\frac{\partial }{\partial x_{j}}=\Delta +\sum_{j=1}^{n}x_{j}\frac{%
\partial }{\partial x_{j}}=\Delta +A 
\]

\noindent where $(\frac{\partial }{\partial x_{j}})^{\ast }$ denotes the
adjoint on $L^{2}(\mathbb{R}^{n},pdx)$ and $A$ is the generator of dilations
on $\mathbb{R}^{n}.$ On $L^{q}(\mathbb{R}^{n},pdx),1<q<\infty ,$ the
spectrum of $N_{0}$ is $\mathbb{N},$ and the Hermite polynomials on $\mathbb{%
R}^{n}$ form an orthogonal basis of eigenvectors of $e^{-tN_{0}}$ in $L^{2}(%
\mathbb{R}^{n},pdx).$

The generator $N_{0}$ has a fruitful generalization in (commutative or non
commutative) analysis on deformed or $q$-Fock spaces, namely the number
operator $N,$ i.e. the second differential quantization of identity. A
substitute of Mehler formula holds and $(e^{-tN})_{t>0}$ is the compression
of a one parameter group of unitary dilations, see e.g. \cite{LP2}.

Our motivation in this paper is to exploit Mehler formula in another
direction: in the setting of stratified groups Mehler formula still defines
a semi-group $(e^{-tN})_{t>0}$ and we study which properties of the
classical O-U semi-group remain valid. We also hope that this semi-group
might throw some light on properties of the heat density $p$.

\textbf{Results and organization of the paper}

In section 2 we recall some properties of the self-adjoint semi-group on $%
L^{2}(pd\gamma )$ whose generator is $-\nabla ^{\ast }\nabla
=-\sum_{i=1}^{n}X_{i}^{\ast }X_{i},$ $X_{i}^{\ast }$ being the formal
adjoint of $X_{i}$ with respect to $L^{2}(pd\gamma ).$ We give in passing a
simple proof of the known Poincar\'{e} inequality in $L^{2}(pd\gamma )$.

In the main section 3 we consider another generalization, the Mehler
semi-group, which is defined for $t\geq 0$ by (theorem \ref{P2})

\[
T_{t}(f)(\gamma )=\int_{G}f(\delta _{e^{-t}}\gamma \delta _{\sqrt{1-e^{-2t}}%
}g)p(g)dg=e^{-tN}(f)(\gamma ). 
\]

\noindent Some properties are described in 3.2, in particular $pd\gamma $ is
an invariant measure. This semi-group is not selfadjoint on $L^{2}(pd\gamma
),$ but formally the real part of its generator $-N$ is $-\nabla ^{\ast
}\nabla $ and $N=L+A$ where $A$ is the generator of the group ($\delta
_{e^{t}})_{t\in \mathbb{R}}$ of dilations$,$ studied in 3.3.

We show in 3.4 that every $T_{t},t>0,$ is compact on $L^{q}(pd\gamma
),1<q<\infty ,$ (proposition \ref{P3}), with common spectrum $e^{-t\mathbb{N}%
}$ on the closed subspace spanned by polynomials (theorem \ref{spectre}),
which coincides with the whole space only if the number of layers of $%
\mathcal{G}$ is $\leq 4$ (proposition \ref{dense}). We describe the
eigenspaces in this case.

In 3.5 we give another description of these eigenspaces if $G$ is step two,
similar to the usual definition of one variable Hermite polynomials by their
generating function.

\noindent\ \ \ \ \ \ \textbf{More} \textbf{notation}

We denote $\mathcal{G=}V_{1}\oplus ..\oplus V_{k}$, where $V_{1},..,V_{k}$
are the layers of the Lie algebra $\mathcal{G}$ of $G,$ $V_{k}=\mathcal{Z}$
being the central layer, so that \cite[p. 5]{FS}

\[
\lbrack V_{j},V_{h}]\subset V_{j+h},\;[V_{1},V_{h}]=V_{h+1},1\leq h<k 
\]

\noindent The homogeneous dimension of $G$ is

\[
Q=\sum_{j=1}^{k}j\dim V_{j}. 
\]

\noindent Generic elements of the layers are denoted respectively by $%
X,Y...,U,$ and respective basis of the layers are denoted by $%
(X_{1},..,X_{n})$, $(Y_{1},..,Y_{m})$, ..., $(U_{1},..,U_{k}).$ Such a basis
is also denoted by $(Z_{j})_{1\leq j\leq N}.$ We denote accordingly

\begin{eqnarray*}
g &=&\exp (\sum x_{i}X_{i}+\sum y_{i}Y_{i}+..+\sum u_{i}U_{i})=\exp
(X+Y+..+U) \\
&=&(x,y,..,u)=\exp (\sum_{j=1}^{N}z_{j}Z_{j})=(z_{j})_{j=1}^{N},
\end{eqnarray*}

\noindent since the mapping $(z_{j})_{j=1}^{N}\rightarrow g$ is a
diffeomorphism: $\mathbb{R}^{N}\rightarrow G.$

We denote by $\mathcal{P}$ the space of polynomials on $G,$ as defined in 
\cite[chapter I-C]{FS} for the fixed basis $(Z_{j})_{j=1}^{N}$: they are
polynomials w.r. to the coordinates $z_{j},1\leq j\leq N.$

The dilation $\delta _{t},t\geq 0,$ are defined on $\mathcal{G}$ and $G$ by

\[
\delta _{t}(X+Y+..+U)=tX+t^{2}Y+..+t^{k}U,\;\;\delta _{t}(\exp Z)=\exp
\delta _{t}(Z),\;Z\in \mathcal{G}. 
\]

\noindent For a function $f$ on $G,$

\[
\delta _{t}(f)=f\circ \delta _{t}. 
\]

\noindent The generator $A$ of the one parameter group ($\delta
_{e^{s}})_{s\in \mathbb{R}}$ of dilations on $G$ satisfies$:$ for $f\in 
\mathcal{S}(G)$ and $s>0$

\begin{equation}
\frac{d}{dt}\mid _{t=1}f\circ \delta _{t}=A(f)=-tt^{A}\frac{d}{dt}%
t^{-A}(f)=-t\delta _{t}\frac{d}{dt}(f\circ \delta _{\frac{1}{t}}).
\label{dil}
\end{equation}

\textbf{Acknowledgment: }We thank W. Hebisch who gave us the idea of the
proof of proposition \ref{dense}.

\section{The semi-group $e^{-t\protect\nabla ^{\ast }\protect\nabla }$ on $%
L^{2}(pdg)$}

This semi-group has already been introduced in \cite{BHT}, under a
probabilistic point of view, in connection with some Markov processes on Lie
groups. \ We use instead an analytic point of view as in \cite{O}. We
consider this semi-group firstly because it is a natural generalization of
the classical O-U semi-group, secondly because its generator $\nabla ^{\ast
}\nabla $ is the real part of the generator $N$ we shall study in part 3,
see theorem \ref{P2}.

\subsection{Definition and some properties}

We consider the (closed) accretive sesquilinear form

\[
a(f,h)=\int_{G}(\nabla f.\nabla h)pdg=\int_{G}\sum_{i=1}^{n}X_{i}f\overline{%
X_{i}h}pdg 
\]

\noindent\ whose (dense) domain in $L^{2}(pdg)$ is the Hilbert space

\[
H^{1}(p)=\{f\in L^{2}(pdg)\mid X_{i}f\in L^{2}(pdg),1\leq i\leq n\} 
\]

\noindent equipped with the norm $\left\Vert f\right\Vert
_{H^{1}(p)}^{2}=\left\Vert f\right\Vert _{L^{2}(p)}^{2}+\left\Vert
\left\vert \nabla f\right\vert \right\Vert _{L^{2}(p)}^{2}$; this form is
continuous on $H^{1}(p)\times H^{1}(p).$

Hence \cite[proposition 1.51, theorem 1.53]{O} it defines an operator, which
we denote by $\nabla ^{\ast }\nabla ,$ such that $-\nabla ^{\ast }\nabla $
is the generator of a strongly continuous semi-group of contractions on $%
L^{2}(pdg);$ moreover this semi-group is holomorphic on the sector $\Sigma _{%
\frac{\pi }{2}}=\{\left\vert \arg z\right\vert <\frac{\pi }{2},z\neq 0\},$
and $e^{-z\nabla ^{\ast }\nabla }$ is a contraction on $L^{2}(pdg)$ for $%
z\in \Sigma _{\frac{\pi }{2}}$. Obviously, on $\mathcal{S}(G),$

\begin{equation}
\nabla ^{\ast }\nabla =\sum_{i=1}^{n}X_{i}^{\ast }X_{i}=L-\sum_{i=1}^{n}%
\frac{X_{i}p}{p}X_{i}=L-B.  \label{500}
\end{equation}

Since $X_{i}$ is a derivation, the chain rule holds$,$ hence $%
X_{i}(f^{+})=(X_{i}f)1_{\{f>0\}}$ by the same proof as for usual derivations
on $\mathbb{R}^{N}$ \cite[proposition 4.4]{O}, and $a(f^{+},f^{-})=0$; since
the form $a$ also preserves real valued functions, the semi-group $%
e^{-t\nabla ^{\ast }\nabla }$ \ is positivity preserving \cite[theorem 2.6]%
{O}. Since $e^{-t\nabla ^{\ast }\nabla }(1)=1,$ the semi-group is thus
contracting on $L^{\infty }(pdg).$ Since moreover $\nabla ^{\ast }\nabla $
is self-adjoint, $e^{-t\nabla ^{\ast }\nabla }$ is measure preserving, i.e.

\[
\int_{G}e^{-t\nabla ^{\ast }\nabla }(f)pdg=\int_{G}fpdg,t>0, 
\]

\noindent so it extends as a contraction semi-group on $L^{1}(pdg)$ hence on 
$L^{q}(pdg),1<q<\infty $ by interpolation.

\subsection{Poincar\'{e} inequality in $L^{2}(pdg)$}

Poincar\'{e} inequality \cite[theorem 4.2]{DM} means that the spectrum of $%
\nabla ^{\ast }\nabla $ on $L^{2}(pdg)$ lies in $\{0\}\cup \lbrack
C^{-1},\infty \lbrack $: there exists $C>0$ such that, for $f\in \mathcal{S}%
(G),$

\begin{equation}
\left\Vert f-\int_{G}fpdg\right\Vert _{L^{2}(pdg)}^{2}\leq
C\int_{G}\left\vert \nabla f\right\vert ^{2}pdg=C\int_{G}f(\nabla ^{\ast
}\nabla f)pdg.  \label{P}
\end{equation}

\noindent (\ref{P}) follows from the inequality (used for $q=2)$ \cite[%
theorem 4.1]{DM}

\begin{equation}
\left\vert \nabla (e^{-tL}f)\right\vert ^{q}\leq C_{q}e^{-tL}(\left\vert
\nabla f\right\vert ^{q}),\;1<q<\infty ,  \label{DM}
\end{equation}

\noindent which B.\ Driver and T. Melcher proved, first for $\mathbb{H}_{1}, 
$ then for nilpotent groups $G$ (see T. Melcher's thesis), using Malliavin
calculus. See also \cite{BHT} for some extensions.

We shall show in proposition \ref{Poincare} that (\ref{DM}) also follows
easily from gaussian estimates of $p$ and $\nabla p$.

Using the explicit formula for the Carnot-Caratheodory distance, H.Q. Li 
\cite[corollary 1.2]{Li}\ obtained (\ref{DM}) for $q=1$, on the
3-dimensional Heisenberg group $G=\mathbb{H}_{1}.$ As well known \cite[th%
\'{e}or\`{e}me 5.4.7]{A}, this implies Log-Sobolev inequality for the
measure $pdg$ on $\mathbb{H}_{1}$ and (\ref{P}). Another proof of this
Log-Sobolev inequality for $\mathbb{H}_{1},$ hence for $\mathbb{H}_{k},$ is
given in \cite[theorem 7.3]{HZ}$.$

\begin{proposition}
\label{Poincare}\cite{DM} Let $G$ be a stratified group. Then (\ref{DM}) and
Poincar\'{e} inequality (\ref{P}) hold true.
\end{proposition}

\noindent Proof: By \cite[theorem 4.2, proposition 2.6, lemma 2.3]{DM} it is
enough to prove (\ref{DM}) for $t=\frac{1}{2}$, at $\gamma =0$. Hence, it is
enough to prove, for an element $X$ of the basis of $V_{1},$ and $f\in 
\mathcal{S}(G),$

\[
\left\vert X(e^{-\frac{1}{2}L}f)(0)\right\vert =\left\vert X(f\ast
p)(0)\right\vert =\left\vert \int_{G}(\widehat{X}f)(g)p(g)dg\right\vert \leq
C_{q,X}\left\Vert \nabla f\right\Vert _{L^{q}(pdg)}; 
\]

\noindent here \cite[p. 22 and proposition 1.29]{FS}

\[
(\widehat{X}f)(g)=\frac{d}{dt}\mid _{t=0}f((\exp tX)g),\;\;\widehat{X}%
=X+\sum_{j>n}Q_{X,j}Z_{j} 
\]

\noindent where $(Z_{j})_{j=1}^{N}$ is a basis of $\mathcal{G}$ respecting
the layers and $Q_{X,j}$ is a polynomial (with homogeneous degree $h-1$ if $%
Z_{j}\in V_{h},2\leq h\leq k).$

\noindent Since $[V_{1},V_{h-1}]=V_{h},2\leq h\leq k,$ we may choose $%
Z_{j}\in V_{h}$ such that $Z_{j}=[Y,A],$ where $Y$ is an element of the
basis of $V_{1}$ and $A\in V_{h-1}.$Then

\[
\left\vert \int_{G}Z_{j}f(g)Q_{X,j}(g)p(g)dg\right\vert \leq \left\vert
\int_{G}Yf\;A(Q_{X,j}p)dg\right\vert +\left\vert
\int_{G}Af\;Y(Q_{X,j}p)dg\right\vert . 
\]

\noindent Iterating for $A\in V_{1}+..+V_{k-1}$ and so on, $\left\vert
\int_{G}(\widehat{X}f)(g)p(g)dg\right\vert $ is finally less than a finite
number (which does not depend on $f)$ of terms $\left\vert
\int_{G}Yf\;Z(Qp)dg\right\vert $ where $Y$ is an element of the basis of $%
V_{1},Z\in \mathcal{G},$ and $Q$ is a polynomial. Each of these terms can be
estimated by

\[
\left\vert \int_{G}Yf\;Z(Qp)dg\right\vert \leq \left\Vert \left\vert \nabla
f\right\vert \right\Vert _{L^{q}(pdg)}(\left\Vert ZQ\right\Vert
_{L^{q^{\prime }}(pdg)}+\left\Vert Q\frac{Zp}{p}\right\Vert _{L^{q^{\prime
}}(pdg)}) 
\]

\noindent where $\frac{1}{q}+\frac{1}{q^{\prime }}=1.$ Then $\left\Vert
ZQ\right\Vert _{L^{q^{\prime }}(pdg)}$ is finite since $ZQ$ is a polynomial
and $p\in \mathcal{S}(G).$ The main point is that $\left\Vert Q\frac{Zp}{p}%
\right\Vert _{L^{q^{\prime }}(pdg)}$ is finite. Indeed, denoting $%
d(g)=d(0,g) $ where $d$ is the Carnot-Caratheodory distance on $G,$ one uses 
\cite[theorem IV.4.2 and Comments on chapter IV]{CSV}: for $0<$ $\varepsilon
<1,$

\begin{equation}
C_{\varepsilon }e^{-\frac{1}{2-2\varepsilon }d^{2}(g)}\leq p(g)\leq
K_{\varepsilon }e^{-\frac{1}{2+2\varepsilon }d^{2}(g)}.  \label{11}
\end{equation}

\noindent and, for $Z\in \mathcal{G},$

\begin{equation}
(Zp)(g)\leq K_{\varepsilon ,Z}\;e^{-\frac{1}{2+2\varepsilon }d^{2}(g)}.
\label{111}
\end{equation}

\noindent\ Hence $Q\frac{Zp}{p}$ lies in $L^{r}(pdg),1\leq r<\infty ,$ which
ends the proof.$\blacksquare $

\section{Definition and properties of the Mehler semi-group}

\subsection{Preliminaries}

\noindent The next proposition extends a classical property of independant
gaussian variables and will imply the semi-group property of our family of
operators.

\begin{proposition}
\label{P1}Let $\gamma ,g$ be independant $G-$valued random variables with
law $pdg.$ Then the r.v.%
\[
\delta _{\cos \theta }\gamma \delta _{\sin \theta }g,\;0\leq \theta \leq 
\frac{\pi }{2} 
\]

\noindent has the same law, i.e. for any bounded borelian function $f$ on $%
G, $

\[
\int_{G^{2}}f(\delta _{\cos \theta }\gamma \delta _{\sin \theta }g)p(\gamma
)p(g)d\gamma dg=\int_{G}f(g)p(g)dg. 
\]

More generally, if $g_{1},...,g_{n}$ are $G$-valued i.i.d r.v. with law $pdg$
and $\sum\limits_{1\leq j\leq n}a_{j}^{2}=1,$ $(a_{j}\geq 0),$ the law of $%
\prod\limits_{j=1}^{j=n}\delta _{a_{j}}g_{j}$ is $pdg.$
\end{proposition}

\noindent Proof: By two changes of variables, denoting $C=\sin \theta \cos
\theta ,$

\begin{eqnarray*}
\int_{G^{2}}f(\delta _{\cos \theta }\gamma \delta _{\sin \theta
}g)p(g)p(\gamma )d\gamma dg &=&\frac{1}{C^{Q}}\int_{G^{2}}f(\gamma ^{\prime
}g^{\prime })p(\delta _{\frac{1}{\cos \theta }}\gamma ^{\prime })p(\delta _{%
\frac{1}{\sin \theta }}g^{\prime })d\gamma ^{\prime }dg^{\prime } \\
&=&\frac{1}{C^{Q}}\int_{G^{2}}f(g)p(\delta _{\frac{1}{\cos \theta }}\gamma
^{\prime })p(\delta _{\frac{1}{\sin \theta }}(\gamma ^{\prime -1}g))d\gamma
^{\prime }dg \\
&=&\int_{G}f(g)(p_{\cos ^{2}\theta }\ast p_{\sin ^{2}\theta })(g)dg \\
&=&\int_{G}f(g)p(g)dg.
\end{eqnarray*}

\noindent The second assertion follows by iteration.

\noindent\ \ \ \ \ \emph{Remark 1}: A central limit theorem for i.i.d
centered random variables with values in a stratified group $G$ and law $\mu 
$ with order 2 moments is proved in \cite[theorem 3.1]{CR}. The density $p$
of the limit law is the kernel at time 1 of a diffusion semi-group whose
generator satisfies (\ref{1}).

\emph{Remark 2: }If $X,Y$ are i.i.d standard gaussian vectors with values in 
$\mathbb{R}^{n}$, the couple $(X\cos \theta +Y\sin \theta ,$ $\frac{d}{%
d\theta }(X\cos \theta +Y\sin \theta ))$ has the same joint law as $(X,Y)$.
This fact\ implies, in the O-U case,\emph{\ }$\ $that $\cos ^{N_{0}}\theta $
is the compression of the isometry $R_{\theta }$ of $L^{2}(\mathbb{R}%
^{n}\times \mathbb{R}^{n},p(x)p(y)dxdy)$ defined by

\[
R_{\theta }(F)(x,y)=F(x\cos \theta +y\sin \theta ,-x\sin \theta +y\cos
\theta ) 
\]

\noindent and $(R_{\theta })_{\theta \in \mathbb{R}}$ is a one parameter
group preserving the measure $p(x)p(y)dxdy.$ This point of view was
exploited\ e.g. in \cite[theorem 2.2]{Q} in order to get a concentration
inequality for the gaussian measure .

In the stratified setting we were not able to exhibit explicit unitary
dilations for the Mehler operators $T_{t}$ defined below.

\subsection{\protect\medskip The Mehler semi-group}

We now define the Mehler semi-group on $L^{q}(G,pdg).$

\begin{theorem}
\label{P2}Let \ $L,$ defined by (\ref{0}), be a subLaplacian on a stratified
group $G,$ and let $p$ be the kernel of $e^{-\frac{L}{2}}$.

a) The family of operators $(T_{t})_{t\geq 0}$ defined on $\mathcal{S(}G%
\mathcal{)}$ by 
\begin{equation}
T_{t}(f)(\gamma )=\int_{G}f(\delta _{e^{-t}}\gamma \delta _{\sqrt{1-e^{-2t}}%
}g)p(g)dg=e^{-\frac{L}{2}(1-e^{-2t})}(f)(\delta _{e^{-t}}\gamma )
\label{1bis}
\end{equation}

is a semi-group whose generator $-N$ is defined on $\mathcal{S(}G\mathcal{)}$
by%
\begin{equation}
N=L+A.  \label{7bis}
\end{equation}

b)The probability measure $pd\gamma $ is invariant by $(T_{t})_{t\geq 0}$
i.e.%
\begin{equation}
\int_{G}T_{t}(f)(\gamma )p(\gamma )d\gamma =\int_{G}f(\gamma )p(\gamma
)d\gamma  \label{3bis}
\end{equation}

and, for $f\in \mathcal{S(}G\mathcal{)},\int_{G}(Nf)pdg=0.$

c) $(T_{t})_{t\geq 0}$ extends as a Markovian semi-group of contractions on $%
L^{q}(G,pd\gamma ),1\leq q\leq \infty $, strongly continuous if $q\neq
\infty .$

d) If $f\in L^{q}(pd\gamma ),1\leq q<\infty ,$ 
\[
\left\Vert T_{t}(f)-\int_{G}fpdg\right\Vert _{L^{q}(pd\gamma )}\rightarrow
_{t\rightarrow \infty }0. 
\]

e) $(T_{t})_{t>0}$ is not self-adjoint on $L^{2}(G,pd\gamma )$ as soon as $G$
is not abelian. Formally \ $\nabla ^{\ast }\nabla $ is the real part of $N,$
i.e., for $f,h\in \mathcal{S(}G\mathcal{)},$%
\[
\left\langle Nf,h\right\rangle _{L^{2}(p)}=\left\langle (\nabla ^{\ast
}\nabla +iC)f,h\right\rangle _{L^{2}(p)} 
\]

where $C$ is a non zero first order differential operator satisfying $%
\left\langle Cf,h\right\rangle =\left\langle f,Ch\right\rangle .$ In
particular, for $f\in \mathcal{S(}G\mathcal{)},$ 
\[
\Re \int_{G}(Nf)fpd\gamma =\int_{G}\left\vert \nabla f\right\vert
^{2}pd\gamma =\int_{G}(\nabla ^{\ast }\nabla f)fpd\gamma . 
\]

If \ moreover $f$ is real valued, the left integral is real.
\end{theorem}

\noindent By the change of notation $e^{-t}=\cos \theta $, $<\theta <\frac{%
\pi }{2}$, (\ref{1bis}) can be rewritten as

\begin{equation}
\cos ^{N}\theta (f)(\gamma )=\int_{G}f(\delta _{\cos \theta }\gamma \delta
_{\sin \theta }g)p(g)dg=\delta _{\cos \theta }\circ e^{-\frac{1}{2}\sin
^{2}\theta L}(f)(\gamma ).  \label{7}
\end{equation}

\noindent Proof: a) Let $\varphi (g^{\prime })=$ $T_{t}(f)(g^{\prime });$ we
compute

\begin{eqnarray*}
T_{s}(\varphi )(\gamma ) &=&\int_{G}\varphi (\delta _{e^{-s}}\gamma \;\delta
_{\sqrt{1-e^{-2s}}}h)p(h)dh \\
&=&\int_{G^{2}}f(\delta _{e^{-t}}[\delta _{e^{-s}}\gamma \delta _{\sqrt{%
1-e^{-2s}}}h]\;\delta _{\sqrt{1-e^{-2t}}}g)p(g)p(h)dgdh \\
&=&\int_{G}f(\delta _{e^{-(t+s)}}\gamma \delta _{\sqrt{1-e^{-2(s+t)}}%
}k)p(k)dk=T_{s+t}(f)(\gamma )
\end{eqnarray*}

\noindent where the third equality comes from proposition \ref{P1} applied
to $(h,g)$.

\noindent By the chain rule applied to (\ref{1bis}),

\[
Nf=-\frac{d}{dt}\mid _{t=0}T_{t}(f)=Lf+A(f). 
\]

\noindent b) Proposition \ref{P1} gives (\ref{3bis}). Differentiating (\ref%
{3bis}) at $t=0$ for $f\in \mathcal{S}(G)$ implies

\[
\int_{G}(Nf)pdg=0. 
\]

\noindent Another proof will be given in Remark 3.

\noindent c) $T_{t}$ is contracting both on $L^{1}(G,pd\gamma ),$ since it
is positivity and measure preserving, and on $L^{\infty }(G,pd\gamma ),$
since it is positivity preserving and $T_{t}(1)=1.$ Hence $T_{t}$ is
contracting on $L^{q}(G,pd\gamma ),1\leq q\leq \infty $ by interpolation.

\noindent\ \ \ \ \ Since $\mathcal{D}(G)$ is norm dense in $L^{q}(G),%
\mathcal{\ }$it is norm dense in $L^{q}(pd\gamma ),1\leq q<\infty :$ indeed,
if $F\in L^{q^{\prime }}(pd\gamma )$ $(\frac{1}{q}+\frac{1}{q^{\prime }}=1)$
and $\int_{G}fFpd\gamma =0$ for every $f\in \mathcal{D}(G),$ then $Fp\in
L^{q^{\prime }}(G)$ hence $Fp=0$ $d\gamma $ a.s.$.$ Writing $e^{-t}=\cos
\theta ,$ one has, for $f\in \mathcal{D}(G),$

\begin{eqnarray*}
\left\Vert T_{t}(f)-f\right\Vert _{L^{q}(pd\gamma )}^{q} &=&\left\Vert
\int_{G}[f(\delta _{\cos \theta }\gamma \delta _{\sin \theta }g)-f(\gamma
)]p(g)dg\right\Vert _{L^{q}(pd\gamma )}^{q} \\
&\leq &\int_{G^{2}}\left\vert f(\delta _{\cos \theta }\gamma \delta _{\sin
\theta }g)-f(\gamma )\right\vert ^{q}p(\gamma )p(g)d\gamma dg,
\end{eqnarray*}

\noindent which converges to 0 as $\theta \rightarrow 0$ by the dominated
convergence theorem. Since $T_{t}$ is contracting, the strong continuity on $%
L^{q}(pd\gamma )$ follows by density.

d) Similarly, if $f$ is bounded and continuous on $G,$

\[
f(\delta _{e^{-t}}\gamma \;\delta _{\sqrt{1-e^{-2t}}}g)\rightarrow
_{t\rightarrow \infty }f(g); 
\]

\noindent by dominated convergence theorem $T_{t}(f)\rightarrow
_{t\rightarrow \infty }\int_{G}f(g)p(g)dg$ \ pointwise and in the norm of $%
L^{q}(pd\gamma ).$ The claim follows by density.

e) By (\ref{7bis}), (\ref{500}) and lemma \ref{20} below, for $f\in \mathcal{%
S}(G),$

\[
(N-\nabla ^{\ast }\nabla )f=A(f)+\sum\limits_{1\leq j\leq n}\frac{X_{j}p}{p}%
X_{j}f=\sum\limits_{1\leq j\leq N}b_{j}Z_{j}f 
\]

\noindent\ where the functions $b_{j}$ are not all zero if $j>n=\dim V_{1}.$
Hence for $h\in \mathcal{S}(G),$

\[
\int_{G}(N-\nabla ^{\ast }\nabla )(f)\overline{h}pdg=-\int_{G}f[\sum%
\limits_{1\leq j\leq N}b_{j}(g)(Z_{j}\overline{h})p+\overline{h}%
Z_{j}(b_{j}p)]dg. 
\]

\noindent By b), the left hand side is zero for $h=1,$ hence $%
\sum\limits_{1\leq j\leq N}Z_{j}(b_{j}p)=0.$ Since $T_{t}$ preserves real
valued functions, so does $N,$ hence

\[
\int_{G}(N-\nabla ^{\ast }\nabla )(f)\overline{h}pdg=-\int_{G}f(N-\nabla
^{\ast }\nabla )(\overline{h})pdg=-\int_{G}f\overline{(N-\nabla ^{\ast
}\nabla )(h)}pdg, 
\]

\noindent which proves ($iC)^{\ast }=-iC,$ where $iC=N-\nabla ^{\ast }\nabla
=A+B$. The remaining assertions are obvious.

\emph{Remark 3}: We now give another instructive proof of $%
\int_{G}(Nf)pdg=0,f\in \mathcal{S}(G),$ hence of (\ref{3bis}). We claim
that, for $f,h\in \mathcal{S}(G)$,

\[
\int_{G}(Nf)hdg=\int_{G}f[L(h)-Qh+\frac{d}{ds}\mid _{s=1}h\circ \delta _{%
\frac{1}{s}}]dg=\int_{G}f(L-QId-A)(h)dg. 
\]

\noindent Indeed, $N=L+A,$ $L$ is formally selfadjoint on $L^{2}(dg)$ and
the claim follows by differentiating at $s=1$ the right hand side of

\[
\int_{G}f(\delta _{s}\gamma )h(\gamma )d\gamma =s^{-Q}\int_{G}f(\gamma
^{\prime })h(\delta _{\frac{1}{s}}\gamma )d\gamma ^{\prime }. 
\]

\noindent By (\ref{dil}) and \cite[lemma 2]{LP}, $p$ may be precisely
defined as the unique solution in $L^{1}(G),$ satisfying $\int_{G}p(g)dg=1,$
of 
\[
(L-QId-A)(p)=Lp-Qp+s\delta _{s}\frac{d}{ds}(p\circ \delta _{\frac{1}{s}%
})=0.\blacksquare 
\]

\emph{Remark 4: } As already mentioned in section 2.2, Log-Sobolev
inequality for $pd\gamma $ is known for $\mathbb{H}_{k}.$ It is equivalent
both to hypercontractivity of $e^{-tN}$ and to hypercontractivity of $%
e^{-t\nabla ^{\ast }\nabla }$ on $\mathbb{H}_{k},$ since $p$ is an invariant
measure for these markovian semigroups and $N,\nabla ^{\ast }\nabla $ are
diffusion operators \cite[theorem 2.8.2]{A}.

\subsection{\protect\smallskip The generator of dilations}

\textbf{\ }We may identify $G$ with a group of finite matrices \cite[theorem
3.6.6]{V}. The derivation formula for an exponential of a matrix valued
function, see e.g. \cite[theorem 69]{H}, applied to a smooth function $Z(s)$%
: $\mathbb{R}\rightarrow \mathcal{G},$ where $\mathcal{G}$ has $k$ layers,
gives

\begin{eqnarray}
\frac{d}{ds}\exp Z(s) &=&\lim_{h\rightarrow 0}\frac{\exp Z(s+h)-\exp Z(s)}{h}
\nonumber \\
&=&\lim_{h\rightarrow 0}\frac{\exp (Z(s)+hZ^{\prime }(s))-\exp Z(s)}{h} 
\nonumber \\
&=&[\exp Z(s)]V(Z(s)),  \label{100}
\end{eqnarray}

\noindent where

\begin{equation}
V(Z(s))=(d\exp )_{Z(s)}(Z^{\prime }(s))=Z^{\prime }(s)+\sum_{l=1}^{k-1}\frac{%
(-1)^{l}}{(l+1)!}(AdZ(s))^{l}(Z^{\prime }(s)).  \label{501}
\end{equation}

\noindent Hence

\[
\exp Z(s+h)=\exp Z(s)\exp h[V(Z(s))+o(1)], 
\]

\noindent which entails for $f\in \mathcal{C}^{\infty }(G)$

\begin{equation}
\frac{d}{ds}f(\exp Z(s))=V(Z(s))(f)(\exp Z(s)).  \label{101}
\end{equation}

\begin{lemma}
\label{20} Let $A$ be the generator of the group of dilations ($\delta
_{e^{t}})_{t\in \mathbb{R}}.$ Then%
\[
A(f)(g)=\sum\limits_{1\leq j\leq N}a_{j}(g)Z_{j}f(g) 
\]

where the functions $a_{j}$ are polynomials w.r. to the coordinates of $g,$
and are not all zero for $j>n=\dim V_{1}$.
\end{lemma}

\noindent Proof: Assume that $\mathcal{G}$ has $k$ layers, $k\geq 2$. Let

\[
\delta _{s}g=\exp (sX+s^{2}Y+..+s^{k}U)=\exp Z(s). 
\]

\noindent By (\ref{101}) $A=V(Z(1)).$ Noting that $Z^{\prime }-Z\in
V_{2}+..+V_{k},$ we get $(AdZ(1))^{l}(Z^{\prime }(1))\in
V_{3}+..+V_{k},l\geq 1.$ So $V(Z(1))-(X+2Y)$ lies in $V_{3}+..+V_{k}$.$%
\blacksquare $

\textbf{Notation: }We denote by $\mathcal{P}_{n}$ the (finite dimensional)
space of homogeneous polynomials on $G$ with homogeneous degree $n,n\in 
\mathbb{N},$ i.e. satisfying

\begin{equation}
\delta _{s}(P)=s^{n}P,\;P\in \mathcal{P}_{n};  \label{vpA}
\end{equation}

\noindent equivalently, $\mathcal{P}_{n}$ is the eigensubspace of $A$ on $%
\mathcal{P}$ associated to $n.$ The finite dimensional subspaces $B_{n}=%
\mathcal{P}_{0}+..+\mathcal{P}_{n}$ are stable under $L$ and dilations,
hence under $e^{-\frac{tL}{2}}$ and $\cos ^{N}\theta $ by (\ref{1bis}),
these operators being naturally extended on $\mathcal{S}^{\prime }(G).$ In
particular $e^{\frac{L}{2}}$ is well defined on $B_{n}$ and is the inverse
of $e^{-\frac{L}{2}},$ which is thus one to one on every $B_{n}$ hence on $%
\mathcal{P=\cup }_{n\geq 0}B_{n}$.

The next lemma is the key for the computation of the spectrum of $\cos
^{N}\theta .$ It will be exploited again in section 3.5.

\begin{lemma}
\label{LA} a) The generator $A$ of dilations on $G$ satisfies $[L,A]=2L$ on $%
\mathcal{C}^{\infty }(G).$

b) $e^{-\frac{L}{2}}\circ \cos ^{N}\theta =\delta _{\cos \theta }e^{-\frac{L%
}{2}}$ on $\mathcal{S}^{\prime }(G).$

c) The set of polynomials $e^{\frac{L}{2}}(\mathcal{P}_{n})$ is a space of
eigenvectors of $\cos ^{N}\theta $ associated to the eigenvalue $\cos
^{n}\theta ,n\geq 0.$
\end{lemma}

\noindent Proof: a) We rewrite (\ref{1}) as

\[
Le^{tA}=e^{2t}e^{tA}L,\;t\in \mathbb{R}, 
\]

\noindent and a) follows by differentiating at $t=0.$

\noindent b) By (\ref{2}), on $\mathcal{S}(G),$ hence on $\mathcal{S}%
^{\prime }(G),$ for $t>0,$

\begin{equation}
e^{-\frac{t^{2}}{2}L}=\delta _{\frac{1}{t}}\circ e^{-\frac{L}{2}}\circ
\delta _{t}.  \label{800}
\end{equation}

\noindent Hence, on $\mathcal{S}^{\prime }(G)),$ by (\ref{1bis}) and (\ref%
{800}) applied to $t=\cos \theta ,$

\[
e^{-\frac{L}{2}}\circ \cos ^{N}\theta =e^{-\frac{L}{2}}\circ \delta _{\cos
\theta }\circ e^{-\frac{\sin ^{2}\theta }{2}L}=\delta _{\cos \theta }\circ
e^{-\frac{L}{2}}. 
\]

\noindent c) Since $e^{-\frac{L}{2}}$ is invertible on $\mathcal{P},$ and $%
\mathcal{P}$ is stable under $cos^{N}\theta ,$ b) implies on $\mathcal{P}$

\[
\cos ^{N}\theta \circ e^{\frac{L}{2}}=e^{\frac{L}{2}}\circ \delta _{\cos
\theta }. 
\]

\noindent Applying this to $\mathcal{P}_{n}$ proves the result.

\noindent

\subsection{Compacity and spectrum of $\cos ^{N}\protect\theta $ on $L^{q}(pd%
\protect\gamma )$}

\begin{proposition}
\label{P3} Let $\cos ^{N}\theta $ be defined by (\ref{7}). Then

a) $\cos ^{N}\theta $ is a Hilbert-Schmidt operator on $L^{2}(pd\gamma ).$

b) $\cos ^{N}\theta $ is compact on $L^{q}(pd\gamma ),1<q<\infty ;$ its non
zero eigenvalues and corresponding eigenspaces are the same on $%
L^{2}(pd\gamma )$ and $L^{q}(pd\gamma ).$ In particular its spectrum $\sigma
(\cos ^{N}\theta )$ does not depend on $q$ and%
\[
\sigma (\cos ^{N}\theta )=(\cos \theta )^{\sigma (N)}\cup \{0\}. 
\]
\end{proposition}

\noindent Actually, $\cos ^{N}\theta $ is a trace class operator on $%
L^{2}(pd\gamma )$ by a) and the semi-group property of $(e^{-tN})_{t>0}.$

\noindent Proof: a) We must show that the kernel of $\cos ^{N}\theta $ lies
in $L^{2}(G\times G,pd\gamma \otimes pdg).$ For fixed $\gamma $ and $\theta
,0<\theta <\frac{\pi }{2},$

\[
\int_{G}f(\delta _{\cos \theta }\gamma \delta _{\sin \theta }g)p(g)dg=\frac{1%
}{\sin ^{Q}\theta }\int_{G}f(z)p(\delta _{\frac{\cos \theta }{\sin \theta }%
}\gamma ^{-1}\delta _{\frac{1}{\sin \theta }}z)dz, 
\]

\noindent so we must prove the convergence of the integral

\[
I(\theta )=\int_{G^{2}}p^{2}(\delta _{\frac{\cos \theta }{\sin \theta }%
}\gamma ^{-1}\delta _{\frac{1}{\sin \theta }}z)\frac{p(\gamma )}{p(z)}%
dzd\gamma . 
\]

\noindent By the gaussian estimates (\ref{11})

\[
\frac{C_{\varepsilon }}{K_{\varepsilon }^{3}}p^{2}(\delta _{\frac{\cos
\theta }{\sin \theta }}\gamma ^{-1}\delta _{\frac{1}{\sin \theta }}z)\frac{%
p(\gamma )}{p(z)}\leq \exp (\frac{d^{2}(z)}{2-2\varepsilon }-\frac{%
d^{2}(\gamma )}{2+2\varepsilon }-\frac{d^{2}(\delta _{\frac{\cos \theta }{%
\sin \theta }}\gamma ^{-1}\delta _{\frac{1}{\sin \theta }}z)}{1+\varepsilon }%
)=\exp \beta . 
\]

\noindent The Carnot distance $d$ satisfies

\[
d(g)\leq d(\gamma ^{-1}g)+d(\gamma )\;and\;d(\delta _{t}g)=td(g). 
\]

\noindent Hence 
\begin{eqnarray*}
(1+\varepsilon )\beta &\leq &\frac{d^{2}(z)}{2(1-\varepsilon )^{2}}-\frac{%
d^{2}(\gamma )}{2}-(\frac{1}{\sin \theta }d(z)-\frac{\cos \theta }{\sin
\theta }d(\gamma ))^{2} \\
&\leq &d^{2}(z)(\frac{1}{2-4\varepsilon }-\frac{1-\cos \theta }{\sin
^{2}\theta })+d^{2}(\gamma )(\frac{\cos \theta -\cos ^{2}\theta }{\sin
^{2}\theta }-\frac{1}{2}).
\end{eqnarray*}

\noindent Since $\frac{1-\cos \theta }{\sin ^{2}\theta }>\frac{1}{2}$ on $]0,%
\frac{\pi }{2}],$ the coefficient of $d^{2}(\gamma )$ is strictly negative,
and so is the coefficient of $d^{2}(z)$ for small enough $\varepsilon >0$.
Hence, for some $c,C>0,$

\[
I(\theta )\leq C\int \int_{G^{2}}e^{-c(d^{2}(z)+d^{2}(\gamma ))}dzd\gamma
=C(\int_{G}e^{-cd^{2}(z)}dz)^{2}. 
\]

\noindent By the left hand side of (\ref{11}), for small $\varepsilon ,$

\[
C_{\varepsilon }\int_{G}e^{-cd^{2}(z)}dz\leq \int_{G}p^{2c(1-\varepsilon
)}(z)dz, 
\]

\noindent and the last integral is finite since $p\in \mathcal{S}(G)$. This
proves a).

b) By interpolation, since $cos^{N}\theta $ is compact on $L^{2}(p(g)dg)$
and bounded on $L^{\infty }(pdg)$ and $L^{1}(pdg),$ it is compact on $%
L^{q}(pdg),1<q<\infty ,$ with the same spectrum and the same eigenspaces
associated to non zero eigenvalues \cite[theorems 1.6.1 and 1.6.2]{D}.

By the compacity on $L^{q}(pdg),$ the set of these eigenvalues is $\{\cos
^{\lambda }\theta \mid \lambda \in \sigma _{q}(N)\}$ where $\sigma _{q}(N)$
denotes the spectrum of $N$ on $L^{q}(pdg)$ \cite[chap. 34.5, theorem 13]{L}%
. Hence $\sigma _{q}(N)=\sigma _{2}(N)$ is discrete and lies in \{$\lambda
\in \mathbb{C}\mid \Re \lambda \geq 0$\} since $cos^{N}\theta $ is
contracting on $L^{2}(pdg)$ (or since $\Re \left\langle
Nf,f\right\rangle \geq 0).\blacksquare $

\begin{theorem}
\label{spectre}Let $G$ be a step $k$ stratified group.

1) If $k\leq 4$

a) the spectrum of $\cos ^{N}\theta $ on $L^{2}(pdg)$ is $\sigma (\cos
^{N}\theta )=(\cos \theta )^{\mathbb{N}}\cup \{0\}$ and $\sigma (N)=\mathbb{N%
}.$

b) the corresponding eigenspaces $E_{n},n\geq 0,$ (which are not pairwise
orthogonal in $L^{2}(pdg)$) are 
\[
E_{n}=e^{\frac{1}{2}L}(\mathcal{P}_{n}). 
\]

2) If $k>4,$ assertions a) b) remain true for the restriction of $\cos ^{%
\mathbb{N}}\theta $ to the closed subspace $L_{\mathcal{P}}^{2}(pdg)$
spanned by polynomials.
\end{theorem}

\noindent If $k=1$ polynomials in $E_{n}$ are the Hermite polynomials with
degree $n$.

\noindent Proof: 1) follows from 2) and proposition \ref{dense} below.

\noindent 2) We first define $E_{n}$ by $E_{n}=e^{\frac{L}{2}}(\mathcal{P}%
_{n}).$ By lemma \ref{LA}, $E_{n}$ lies in the eigenspace of $\cos
^{N}\theta $ associated to the eigenvalue $\cos ^{n}\theta .$ By proposition %
\ref{P3}, $\cos ^{N}\theta $ is compact on $L_{\mathcal{P}}^{2}(pdg).$ The
claim then follows from the following facts:

Let $T:E\rightarrow E$ be a compact operator on an infinite dimensional
Banach space $E;$ let $\Lambda $ be a set of eigenvalues of $T$ and let $%
E_{\lambda },\lambda \in \Lambda ,$ be eigensubspaces whose union is total
in $E.$ Then

a) the spectrum of $T$ is $\Lambda \cup \{0\}$

b) for $\lambda \in \Lambda ,E_{\lambda }$ is the whole eigenspace
associated to $\lambda .$

\noindent Indeed, assume that $T$ has an eigenvalue $\lambda _{0}\notin
\Lambda .$ Then $T-\lambda _{0}I$ \ has a closed range with non zero finite
codimension (see e.g. \cite[chap. 21.1, theorems 3, 4]{L}).\ But this range
contains the linear span of the $E_{\lambda }$'s, $\lambda \in \Lambda ,$
hence is the whole of $E.$ This is a contradiction, which proves a).

Let $\lambda _{0}\in \Lambda ;$ since $E_{\lambda _{0}}$ is stable under $T,$
$T$ acts on the quotient space $E/E_{\lambda _{0}}$ and is still compact.
The $E_{\lambda }$'s, $\lambda \in \Lambda \backslash \{\lambda _{0}\}$ span
a dense subspace of $E/E_{\lambda _{0}}.$ Applying a) to $E/E_{\lambda _{0}}$%
, $\lambda _{0}$ cannot belong to the spectrum of $T$ on the quotient space,
which proves b).$\blacksquare $

The proof of the next proposition is essentially due to W. Hebisch (private
communication).

\begin{proposition}
\label{dense}Let $G$ be a stratified group. Then the polynomials are dense
in $L^{2}(pdg)$ if and only if $G$ is step $k$ with $k\leq 4$.
\end{proposition}

\noindent Proof: 1) We recall that polynomials are dense in $L^{2}(\mathbb{R}%
,e^{-c\left\vert x\right\vert ^{\alpha }}dx)$ if and only if $\alpha \geq 
\frac{1}{2}:$ obviously, this does not depend on $c$ and is equivalent to
the density of polynomials in $L^{2}(\mathbb{R}^{+},e^{-x^{\alpha }}dx).$ If 
$0<\alpha <\frac{1}{2},$ \cite[Part III, problem 153]{PS} produces a non
zero bounded function $g_{\alpha }$ which is orthogonal to polynomials in $%
L^{2}(\mathbb{R}^{+},e^{-\cos (\alpha \pi )x^{\alpha }}dx).$ If $\alpha \geq 
\frac{1}{2},$ the result follows from the trick of \cite[p 197-198]{Ham}.
Indeed, if $\psi \in L^{2}(\mathbb{R}^{+},e^{-x^{\alpha }}dx)$ and $\alpha
\geq \frac{1}{2}$, the function

\[
F(z)=\int_{\mathbb{R}^{+}}\psi (x)e^{\sqrt{x}z}e^{-x^{\alpha }}dx=\int_{%
\mathbb{R}^{+}}\psi (y^{2})e^{yz}e^{-y^{2\alpha }}ydy 
\]

\noindent is bounded and holomorphic on $\{\Re z<\beta \}$ for some $%
\beta >0,$ by Cauchy-Schwarz inequality.

\noindent Expanding $z\rightarrow e^{\sqrt{x}z}$ in power series, one gets $%
F(-z)=-F(z)$ if $\psi $ is orthogonal to polynomials in $L^{2}(\mathbb{R}%
^{+},e^{-x^{\alpha }}dx)$. Thus $F$ extends as a bounded entire function,
which must be zero by Liouville theorem since $F(0)=0$. Hence the Fourier
transform of $y\rightarrow \psi (y^{2})e^{-y^{2\alpha }}y$ is zero, i.e. $%
\psi =0$ a.s..

2) We identify $g=\exp Z\in G$ with the coordinates $(x,y,..,w)$ of $Z$ \
w.r. to a basis respecting the layers and denote

\[
\eta (g)=\sum_{i\leq l}\left\vert x_{i}\right\vert ^{2}+\sum_{i\leq
m}\left\vert y_{i}\right\vert ^{\frac{2}{2}}+...+\sum_{i\leq r}\left\vert
w_{i}\right\vert ^{\frac{2}{k}}. 
\]

\noindent Obviously $\eta (\delta _{s}g)=s^{2}\eta (g),$ in particular $\eta
(g)=d^{2}(g)\eta (\delta _{\frac{1}{d(g)}}g),$ $d$ denoting the Carnot
distance. Since $\eta $ is strictly positive and bounded on the $d$-unit
sphere of $G,$ there exist constants $c^{\prime },C^{\prime }>0$ such that

\[
c^{\prime }\eta (g)\leq d^{2}(g)\leq C^{\prime }\eta (g). 
\]

\noindent By (\ref{11}) there exist constants $c,C>0$ such that the
following embeddings

\[
L^{2}(e^{-C\eta (g)}dg)\rightarrow L^{2}(pdg)\rightarrow L^{2}(e^{-c\eta
(g)}dg) 
\]

\noindent are continuous, with dense ranges since $\mathcal{D}(G)$ is dense
in the three spaces.

3) The algebraic tensor product

\[
\mathcal{E}=\otimes _{i\leq l}L^{2}(e^{-Cx_{i}^{2}}dx_{i})\otimes ...\otimes
_{i\leq p}L^{2}(e^{-C\left\vert w_{i}\right\vert ^{\frac{2}{k}}}dw_{i}), 
\]%
is dense in $L^{2}(e^{-C\eta (g)}dg).$ For $k\leq 4,$ one variable
polynomials are dense in every factor of $\mathcal{E}$ by step 1), hence
polynomials are dense in $L^{2}(e^{-C\eta (g)}dg)$ and in $L^{2}(pdg).$

Let $k\geq 5.$ By 1) there exists a non zero function $g\in
L^{2}(e^{-c\left\vert w_{r}\right\vert ^{\frac{2}{k}}}dw_{r})$ which is
orthogonal to polynomials w.r. to $w_{r}$. Then $1\otimes ....\otimes
1\otimes g\in L^{2}(e^{-c\eta (g)}dg)$ is orthogonal to all polynomials, so
polynomials are neither dense in $L^{2}(e^{-c\eta (g)}dg),$ nor in $%
L^{2}(pdg).\blacksquare $

\subsection{Generating functions of polynomial eigenvectors of $N$}

The usual Hermite polynomials on $\mathbb{R},$ denoted by $H_{n},n\in 
\mathbb{N},$ are the eigenvectors of the Ornstein-Uhlenbeck operator $N_{0},$
and have the generating function

\[
e^{ixt+\frac{1}{2}t^{2}}=\sum_{n\geq 0}\frac{(it)^{n}}{n!}H_{n}(x)=e^{\frac{1%
}{2}\Delta }(e^{ixt})=e^{\frac{1}{2}\Delta }\circ \delta _{t}(e^{ix}), 
\]

\noindent noting that $x\rightarrow e^{ix}$ is a bounded eigenvector of $%
\Delta .$ In particular

\[
i^{n}H_{n}(x)=\frac{d^{n}}{dt^{n}}\mid _{t=0}e^{\frac{1}{2}\Delta }\circ
\delta _{t}(e^{ix}). 
\]%
We shall verify (proposition \ref{P5}) that a similar formula gives
polynomial eigenvectors of $N.$ When $G$ is step two, these vectors are
total in $L^{q}(pdg),1\leq q<\infty ,$ see theorem \ref{9} below. More
precisely we give in 3.5.1 a technical lemma producing eigenvectors of $N$
out of eigenvectors of $L.$ In 3.5.3 we use this lemma when $\varphi $ is
both an eigenvector of $L$ and a coefficient function of a representation of 
$G$ (proposition \ref{P5}). We shall first gather in 3.5.2 well known facts
about these functions.

\subsubsection{Candidates for generating functions of eigenvectors of $N$}

In the next lemma \ref{GF} we state technical assumptions ensuring the
validity of the computation of some eigenvectors of $N.$ Using lemma \ref{LA}
b), the point is to define $"e^{\frac{L}{2}}\varphi "$ for suitable
functions $\varphi :$ in lemma \ref{LA} c), we choose $\varphi \in \mathcal{P%
},$ here we choose eigenvectors of $L.$

\begin{lemma}
\label{GF}\noindent Let $G$ be a stratified group and let $\varphi \in $ $%
\mathcal{S}^{\prime }(G)\cap \mathcal{C}^{\infty }(G)$ be an eigenvector of $%
L$ such that $L\varphi =\lambda \varphi .$ We assume that, for $n\geq 1,$

(i) $\frac{d^{n}}{dt^{n}}\mid _{t=0}\int_{G}\delta _{t}(\varphi )(\gamma
g^{-1})p(g)dg=\int_{G}\frac{d^{n}}{dt^{n}}\mid _{t=0}\delta _{t}(\varphi
)(\gamma g^{-1})p(g)dg$

(ii) $\frac{d^{n}}{dt^{n}}\mid _{t=0}\delta _{t}(\varphi )$ is a polynomial
on $G.$

Let 
\[
\ f_{t}=e^{\frac{t^{2}\lambda }{2}}\delta _{t}(\varphi ),\;t>0;\;h_{n}=\frac{%
d^{n}}{dt^{n}}\mid _{t=0}f_{t}. 
\]

Then $h_{n}$ is a polynomial on $G$ and 
\[
\cos ^{N}\theta (h_{n})=\cos ^{n}\theta \;h_{n}. 
\]
\end{lemma}

Proof: Since $\varphi \in $ $\mathcal{C}^{\infty }(G)$, $t\rightarrow f_{t}$
is $\mathcal{C}^{\infty }$ on $\mathbb{R}^{+}.$ By (\ref{1}) $L\circ \delta
_{t}(\varphi )=t^{2}\lambda \delta _{t}(\varphi )$, so that $\delta
_{t}(\varphi )=e^{-\frac{L}{2}}f_{t}.$ By lemma \ref{LA} b)

\begin{equation}
e^{-\frac{L}{2}}\cos ^{N}\theta (f_{t})=\delta _{\cos \theta }e^{-\frac{L}{2}%
}f_{t}=\delta _{\cos \theta }\delta _{t}(\varphi )=\delta _{t\cos \theta
}(\varphi )=e^{-\frac{L}{2}}f_{t\cos \theta }.  \label{213}
\end{equation}

\noindent We claim that

\begin{equation}
\frac{d^{n}}{dt^{n}}\mid _{t=0}e^{-\frac{L}{2}}\cos ^{N}\theta (f_{t})=e^{-%
\frac{L}{2}}\cos ^{N}\theta (\frac{d^{n}}{dt^{n}}\mid _{t=0}f_{t})=e^{-\frac{%
L}{2}}\cos ^{N}\theta (h_{n}).  \label{214}
\end{equation}

\noindent In particular, applying (\ref{214}) with $\theta =0,\frac{d^{n}}{%
dt^{n}}\mid _{t=0}e^{-\frac{L}{2}}(f_{t})=e^{-\frac{L}{2}}(h_{n}).$

\noindent Hence, by (\ref{214}) and (\ref{213}),

\begin{equation}
e^{-\frac{L}{2}}\cos ^{N}\theta \;(h_{n})=\frac{d^{n}}{dt^{n}}\mid _{t=0}e^{-%
\frac{L}{2}}f_{t\cos \theta }=e^{-\frac{L}{2}}\cos ^{n}\theta \;h_{n}.
\label{216}
\end{equation}

By Leibnitz rule, it is enough to prove the claim for $\delta _{t}(\varphi )$
instead of $f_{t}.$ By lemma \ref{LA} b) we may replace $e^{-\frac{L}{2}%
}\cos ^{N}\theta $ in the claim by $\delta _{\cos \theta }e^{-\frac{L}{2}}.$
The claim now follows from assumption (i).

By Leibnitz rule and assumption (ii), $h_{n}$ is a polynomial. So is $\cos
^{N}\theta (h_{n})$ and the result follows from (\ref{216}) since $e^{-\frac{%
L}{2}}$ is one to one on $\mathcal{P}.$

\emph{Remark 5: }$\varphi $ and $\varphi \circ \delta _{\beta },\beta >0,$
give colinear $h_{n}$'s, since 
\[
\frac{d^{n}}{dt^{n}}\mid _{t=0}e^{\frac{1}{2}t^{2}\beta ^{2}\lambda }\delta
_{t\beta }(\varphi )=\beta ^{n}\frac{d^{n}}{dt^{n}}\mid _{t=0}e^{\frac{1}{2}%
t^{2}\lambda }\delta _{t}(\varphi )=\beta ^{n}h_{n}. 
\]

\subsubsection{A total set of eigenvectors of $L$ in $L^{q}(pdg),1\leq
q<\infty .$}

Let $\Pi :G\rightarrow B(L^{2}(\mathbb{R}^{k},d\xi ))$ be a non trivial
unitary irreducible representation of $G.$ By definition, $F\in L^{2}(%
\mathbb{R}^{k})$ is a $\mathcal{C}^{\infty }$ vector for $\Pi $ if the
vector valued function: $g\rightarrow \Pi (g)(F)$ is $\mathcal{C}^{\infty }$
on $G.$ We still denote by $\Pi $ the associated differential
representation, defined for a $\mathcal{C}^{\infty }$ vector $F$ and $X\in 
\mathcal{G}$ by

\begin{equation}
X\Pi (g)(F)=\frac{d}{dt}\mid _{t=0}\Pi (g\exp tX)(F)=\Pi (g)\Pi
(X)(F),\;g\in G,  \label{50}
\end{equation}

\noindent and $\Pi (X^{m})=\Pi (X)^{m},$ see e.g. \cite[p.227]{CG}; by
definition, $\Pi (X^{m})(F)$ still lies in $L^{2}(\mathbb{R}^{k})$ and is
still a $\mathcal{C}^{\infty }$ vector for $\Pi .$

\noindent $\Pi $ extends as a representation of the convolution algebra $%
M(G) $ by

\[
\Pi (\mu )=\int\limits_{G}\Pi (g)d\mu (g). 
\]%
In particular\ ($\Pi (p_{t}dg))_{t\geq 0}$ is a semigroup of operators on $%
L^{2}(\mathbb{R}^{k}),$ whose generator is $-\Pi (L).$ Indeed, for a $%
\mathcal{C}^{\infty }$ vector $F$, by (\ref{50}),

\begin{eqnarray*}
-\frac{d}{dt}\int_{G}\Pi (g)(F)p_{t}(g)dg &=&\int_{G}\Pi
(g)(F)(Lp_{t})(g)dg=\int_{G}L\circ \Pi (g)(F)p_{t}(g)dg \\
&=&\int_{G}\Pi (g)\circ \Pi (L)(F)p_{t}(g)dg\rightarrow _{t\rightarrow
0^{+}}\Pi (L)(F).
\end{eqnarray*}

\noindent Since $p\in \mathcal{S}(G)$, $\Pi (pdg)=e^{-\frac{1}{2}\Pi (L)}$
is a trace class operator \cite[theorem 4.2.1]{CG}; in particular its non
zero eigenvalues are \{$e^{-\frac{1}{2}\lambda },\lambda \in \sigma _{2}(\Pi
(L))\},$ where $\lambda $ runs through the eigenvalues of $\Pi (L)$ on $%
L^{2}(\mathbb{R}^{k})$. Moreover, for $F\in L^{2}(\mathbb{R}^{k}),$ the
function $\Pi (pdg)(F)$ is a $\mathcal{C}^{\infty }$ vector for $\Pi $ \cite[%
theorem A.2.7 p. 241]{CG}.

Let $\mathcal{U}$ be a set of non trivial unitary irreducible
representations of $G$ whose equivalence classes support the Plancherel
measure for $G.$ By Kirillov theory, there exists an integer $k,$ which does
not depend on $\Pi \in \mathcal{U}$ , such that $\Pi :G\rightarrow B(L^{2}(%
\mathbb{R}^{k}))$, see more details in 3.5.4 below.

\begin{proposition}
\label{L5}Let $G$ be a stratified group and let $\mathcal{F}$ be the set of
coefficient functions 
\[
\mathcal{F=}\{\varphi ^{\Pi ,\mu ,\mu ^{\prime }}=\left\langle \Pi
(.)(F_{\mu }),F_{\mu ^{\prime }}\right\rangle \mid \Pi \in \mathcal{U}%
,F_{\mu },F_{\mu ^{\prime }}\in \mathcal{B}_{\Pi }\}\subset L^{\infty }(dg)
\]

where $\mathcal{B}_{\Pi }$ is an orthogonal basis of $L^{2}(\mathbb{R}^{k})$
chosen among eigenvectors of $e^{-\frac{1}{2}\Pi (L)}.$ Then $\mathcal{F},$
which lies in $\mathcal{C}^{\infty }(G),$ is a set of eigenvectors of $L$
which is total in $L^{q}(p(g)dg),1\leq q<\infty .$

For fixed $\ \Pi ,\mu $ the functions $\{\varphi ^{\Pi ,\mu ,\mu ^{\prime
}}\mid F_{\mu ^{\prime }}\in \mathcal{B}_{\Pi }\}$ are independent and
belong to the same eigenspace of $L.$
\end{proposition}

\noindent Proof: a) For every non trivial unitary irreducible representation 
$\Pi $ of $G,$ since $\Pi (pdg)(F_{\mu })=e^{-\frac{1}{2}\Pi (L)}(F_{\mu
})=e^{-\frac{1}{2}\lambda _{\mu }}F_{\mu },$ $F_{\mu }$ is a $\mathcal{C}%
^{\infty }$ vector for $\Pi ,$ hence $\varphi ^{\Pi ,\mu ,\mu ^{\prime }}\in 
\mathcal{C}^{\infty }(G)$; $\varphi ^{\Pi ,\mu ,\mu ^{\prime }}$ is an
eigenvector of $L$ with eigenvalue $\lambda _{\mu }$ by (\ref{50}).

\noindent Since $\Pi $ is irreducible, the closed invariant subspace

\[
\{F\in L^{2}(\mathbb{R}^{k})\mid \forall g\in G\;\left\langle \Pi (g)(F_{\mu
}),F\right\rangle =0\} 
\]

\noindent is reduced to $\{0\},$ which implies the independence of the $%
\varphi ^{\Pi ,\mu ,\mu ^{\prime }}$'s. (In the Heisenberg case, see \cite[%
p. 19, 51]{T}).

b) Let $\psi \in L^{q^{\prime }}(pdg),\frac{1}{q}+\frac{1}{q^{\prime }}=1,$
be orthogonal to $\mathcal{F}$, i.e. for $\Pi \in $ $\mathcal{U},$ 
\[
0=\int_{G}\left\langle \Pi (g)(F_{\mu }),F_{\mu ^{\prime }}\right\rangle
\psi (g)p(g)dg=\left\langle (\int_{G}\Pi (g)\psi (g)p(g)dg)(F_{\mu }),F_{\mu
^{\prime }}\right\rangle . 
\]

\noindent Equivalently $\Pi (\psi p)=\widehat{\psi p}(\Pi )=0$ for $\Pi \in 
\mathcal{U}.$ Then Plancherel formula for $G$ (see e.g. \cite[theorem 4.3.10]%
{CG}) implies that $\psi p=0$ $dg$ a.s.. Indeed, this is clear if $\psi p\in
L^{2}(dg),$ in particular if $q^{\prime }\geq 2.$ In general, $\psi p\in
L^{1}(dg),\left\Vert (\psi p)\ast p_{t}-\psi p\right\Vert
_{L^{1}(dg)}\rightarrow _{t\rightarrow 0}0$ and $(\psi p)\ast p_{t}\in
L^{2}(dg);$ moreover $(\psi p)\ast p_{t}=0$ a.s. since, for every $\Pi \in 
\mathcal{U},$

\[
\Pi ((\psi p)\ast p_{t})=\Pi (\psi p)\Pi (p_{t})=0.\blacksquare 
\]

\subsubsection{Polynomial eigenvectors of N built from coefficients of
representations}

We now consider the functions $e^{\frac{1}{2}t^{2}\lambda _{\mu }}\varphi
^{\Pi ,\mu ,\mu ^{\prime }}\circ \delta _{t}$ as generating functions of
polynomial eigenvectors of $N.$

\begin{proposition}
\label{P5} Let $\varphi ^{\Pi ,\mu ,\mu ^{\prime }}=\left\langle \Pi
(.)(F_{\mu }),F_{\mu ^{\prime }}\right\rangle \in \mathcal{F}$ be as in
proposition \ref{L5}. For $n\geq 1,$ let%
\[
h_{n}^{\Pi ,\mu ,\mu ^{\prime }}=\frac{d^{n}}{dt^{n}}\mid _{t=0}e^{\frac{1}{2%
}t^{2}\lambda _{\mu }}\varphi ^{\Pi ,\mu ,\mu ^{\prime }}\circ \delta _{t}. 
\]

Then $h_{n}^{\Pi ,\mu ,\mu ^{\prime }}$ is a polynomial eigenvector of $\cos
^{N}\theta $ with eigenvalue $\cos ^{n}\theta $.
\end{proposition}

\noindent Proof: By proposition \ref{L5} and lemma \ref{GF}, it is enough to
prove assumptions (i) and (ii) in lemma \ref{GF}. We claim the existence of
a polynomial $\psi _{n},$ $n\geq 1,$ which does not depend on $t,$ such
that, for $0\leq t\leq 1$ and $n\geq 0,$%
\[
\;\;\left\vert \frac{d^{n}}{dt^{n}}\;\varphi ^{\Pi ,\mu ,\mu ^{\prime
}}\circ \delta _{t}\right\vert \leq \psi _{n}. 
\]

\noindent\ Since $g\rightarrow \psi _{n}(\gamma g^{-1})$ is still a
polynomial, it lies in $L^{1}(pdg)$ for every $\gamma \in G,$ and this will
prove assumption (i). We now verify the claim.

\textbf{Case 1: }The computation of derivatives being easier if $G$ is step
two, we first consider this setting.

\noindent By Schur lemma, the restriction of $\Pi $ to the center $\exp 
\mathcal{Z}$ of $G$ is given by a character $u\rightarrow e^{i\left\langle
l,u\right\rangle }$ where $l$ is some linear form on $\mathcal{Z},$ see e.g. 
\cite[p. 184]{CG}. If $g=(x,u)$ and $X=\sum_{j=1}^{n}x_{j}X_{j}\in V_{1},$

\[
\varphi ^{\Pi ,\mu ,\mu ^{\prime }}(\delta _{t}g)=e^{it^{2}\left\langle
l,u\right\rangle }\left\langle \Pi (\exp tX)(F_{\mu }),F_{\mu ^{\prime
}}\right\rangle =e^{it^{2}\left\langle l,u\right\rangle }\Phi _{t}^{\Pi ,\mu
,\mu ^{\prime }}(x) 
\]

\noindent and, by (\ref{50}),

\begin{equation}
\frac{d^{m}}{dt^{m}}\Phi _{t}^{\Pi ,\mu ,\mu ^{\prime }}(x)=\left\langle \Pi
(\exp tX)\Pi (X)^{m}(F_{\mu }),F_{\mu ^{\prime }}\right\rangle .  \label{215}
\end{equation}

\noindent Since $\Pi (X)^{m}(F_{\mu })$ lies in $L^{2}(\mathbb{R}^{k},d\xi
),\left\langle \Pi (X)^{m}(F_{\mu }),F_{\mu ^{\prime }}\right\rangle $ and $%
\left\Vert \Pi (X)^{m}(F_{\mu })\right\Vert _{L^{2}(d\xi )}$ are polynomials
w.r. to $x,$ $\frac{d^{m}}{dt^{m}}\mid _{\alpha =0}\delta _{t}(\varphi ^{\Pi
,\mu ,\mu ^{\prime }})$ is a polynomial w.r. to $x,u,$ and $\left\vert \frac{%
d^{m}}{dt^{m}}e^{it^{2}\left\langle l,u\right\rangle }\Phi _{t}^{\Pi ,\mu
,\mu ^{\prime }}(x)\right\vert $ is, for $0\leq t\leq 1,$ less than a
polynomial $\psi _{n}$ which does not depend on $t.$ This proves (i) and
(ii) in this case.

\textbf{General case: } As in (\ref{100}) and (\ref{501}), for $g=\exp
Z=\exp (X+Y+..+U)$ and $t>0,$ since $V(\Pi (\delta _{t}Z))=\Pi (V(\delta
_{t}Z)),$

\[
\frac{d}{dt}\varphi ^{\Pi ,\mu ,\mu ^{\prime }}(\delta _{t}g)=\frac{d}{dt}%
\left\langle \exp \Pi (\delta _{t}Z)(F_{\mu }),F_{\mu ^{\prime
}}\right\rangle =\left\langle \Pi (V(\delta _{t}Z))(F_{\mu }),\exp -\Pi
(\delta _{t}Z)(F_{\mu ^{\prime }})\right\rangle . 
\]

\noindent At $t=0$ this reduces to the polynomial $\left\langle \Pi
(X)(F_{\mu }),F_{\mu ^{\prime }}\right\rangle $. Since $\Pi (V(\delta _{t}Z)$
has polynomial coefficients w.r. to $t$ and the coordinates of $g,$ so does $%
\left\Vert \Pi (V(\delta _{t}Z))(F_{\mu })\right\Vert _{L^{2}(d\xi )}.$
Hence there is a polynomial $\psi _{1}$ w.r. to the coordinates of $g$ such
that

\[
sup_{0\leq t\leq 1}\left\Vert \Pi (V(\delta _{t}Z))(F_{\mu })\right\Vert
_{L^{2}(d\xi )}\leq \psi _{1}. 
\]

\noindent This proves the claim for $n=1$. Clearly this can be iterated for
upper derivatives, which proves (i) and (ii).$\blacksquare $

\subsubsection{The step two setting: generalized Hermite polynomials}

In this case, the key facts are the extension of the explicit functions $%
\varphi ^{\Pi ,\mu ,\mu ^{\prime }}\in \mathcal{F}$ as entire functions on
the complexification of $G$ and the explicit expression of $p.$ Theorem \ref%
{9} gives another proof of theorem \ref{spectre} a) in this setting, with
another description of the eigenspaces of $N$ by generating functions.

\begin{theorem}
\label{9}Let $G$ be a step two stratified group. Then

a) every $\varphi ^{\Pi ,\mu ,\mu ^{\prime }}\in \mathcal{F}$ lies in the
closed subspace of $L^{q}(pdg),1\leq q<\infty ,$ spanned by constants and
the polynomials $\{h_{n}^{\Pi ,\mu ,\mu ^{\prime }},n\geq 1\}$ defined in
proposition \ref{P5}$.$

b) The set of generalized Hermite polynomials%
\[
\cup _{\varphi ^{\Pi ,\mu ,\mu ^{\prime }}\in \mathcal{F}}\{h_{n}^{\Pi ,\mu
,\mu ^{\prime }},n\geq 1\} 
\]

together with the constants is a set of eigenvectors of $N$ which is total
in $L^{q}(pdg),1\leq q<\infty .$

c) For $\ $fixed $n\geq 1,$ $\cup _{\varphi ^{\Pi ,\mu ,\mu ^{\prime }}\in 
\mathcal{F}}\{h_{n}^{\Pi ,\mu ,\mu ^{\prime }}\}$ spans the eigenspace of $N$
associated to $n$ in $L^{q}(pdg),1<q<\infty .$
\end{theorem}

In contrast, if $G$ has more than 4 layers, assertion b) is false by
proposition \ref{dense}, hence a) is false for some $\varphi ^{\Pi ,\mu ,\mu
^{\prime }}\in \mathcal{F}$, by proposition \ref{L5}$.$ If $G$ has 3 or 4
layers, we do not know if the conclusions of theorem \ref{9} hold true.

\noindent Proof of theorem \ref{9}: a) implies b) by propositions \ref{L5}
and \ref{P5}.

b) implies c) as recalled in the proof of theorem \ref{spectre}.

\noindent\ \ \ \ \ \ a) The proof is given in three steps. In step 1 we
state two standard sufficient conditions ensuring statement a); in step 2 we
verify these conditions when $G$ is a Heisenberg group; in step 3 we show
how the general step 2 case mimicks the Heisenberg \ case.

\noindent\ \ \ \ \ \textbf{Step 1:} \ Let $\varphi ^{\Pi ,\mu ,\mu ^{\prime
}}\in \mathcal{F}$ and assume that

(i) for every $g\in G,$ the function $t\rightarrow \varphi ^{\Pi ,\mu ,\mu
^{\prime }}(\delta _{t}g)$ extends as a holomorphic function $z\rightarrow
\varphi _{z}^{\Pi ,\mu ,\mu ^{\prime }}(g)$ on $\mathbb{C}.$

(ii) for some connected neighborhood $\Omega $ of the real axis, for every
compact $K\subset \Omega ,\ $there exists $w_{K}\in L^{q}(pdg),$ $1\leq
q<\infty ,$ such that 
\[
\left\vert \varphi _{z}^{\Pi ,\mu ,\mu ^{\prime }}\right\vert \leq
w_{K},\;z\in K. 
\]

\noindent We claim that $\varphi ^{\Pi ,\mu ,\mu ^{\prime }}=\varphi $ then
lies in the closed subspace of $L^{q}(pdg)$ spanned by $h_{n}^{\Pi ,\mu ,\mu
^{\prime }},n\geq 1.$ Indeed, let $\psi \in L^{q^{\prime }}(pdg),\frac{1}{q}+%
\frac{1}{q^{\prime }}=1,$ and let

\[
m(t)=\int_{G}\varphi (\delta _{t}g)\psi (g)p(g)dg. 
\]

\noindent By the assumptions, $m$ extends as a holomorphic function on $%
\Omega $ and

\[
\frac{d^{n}}{dz^{n}}m=\int_{G}(\frac{d^{m}}{dz^{m}}\varphi _{z})\psi
pdg,\;m\geq 0. 
\]

\noindent By proposition \ref{L5}, $L(\varphi )=\lambda \varphi $ for some $%
\lambda =\lambda _{\mu }.$ Hence $t\rightarrow e^{\frac{1}{2}t^{2}\lambda
}m(t)$ also extends as a holomorphic function on $\Omega $ \ and

\[
\frac{d^{n}}{dz^{n}}\mid _{z=0}e^{\frac{1}{2}z^{2}\lambda }m=\int_{G}[\frac{%
d^{n}}{dz^{n}}\mid _{z=0}e^{\frac{1}{2}z^{2}\lambda }\varphi _{z}]\psi
pdg=\int_{G}h_{n}^{\Pi ,\mu ,\mu ^{\prime }}\psi pdg,\;n\geq 0. 
\]

\noindent If $\psi $ is orthogonal to $\{h_{n}^{\Pi ,\mu ,\mu ^{\prime
}},n\geq 0\},$ these derivatives are zero, hence $e^{\frac{1}{2}z^{2}\lambda
}m$ \ is zero on $\Omega .$ In particular $m(1)=0,$ i.e. $\psi $ is
orthogonal to $\varphi ,$ which proves the claim.

\textbf{Step 2: The Heisenberg groups} $\mathbb{H}_{k}$

\noindent A basis of the first layer of the Lie algebra is $%
X_{1},Y_{1},..,X_{k},Y_{k}$ where $[X_{j},Y_{j}]=-4U,$ $U$ spans the center,
and the other commutators are zero. By the Campbell-Hausdorff formula,

\[
g=\exp (\sum_{j=1}^{k}x_{j}X_{j}+y_{j}Y_{j}+uU)=\exp uU\
\prod\limits_{j=1}^{k}\exp (-2x_{j}y_{j}U)\exp y_{j}Y_{j}\exp x_{j}X_{j}. 
\]

\noindent We first consider the Schr\"{o}dinger (unitary irreducible)
representation $\Pi _{S}:\mathbb{H}_{k}\rightarrow B(L^{2}(\mathbb{R}^{k})),$
defined on the Lie algebra by

\[
\Pi _{S}(X_{j})=\frac{\partial }{\partial \xi _{j}},\;\Pi _{S}(Y_{j})=i\xi
_{j},\;\Pi _{S}(U)=-\frac{1}{4}[\frac{\partial }{\partial \xi _{j}},i\xi
_{j}]=-\frac{i}{4}I. 
\]%
For $F\in L^{2}(\mathbb{R}^{k}),$ this implies

\begin{equation}
\Pi _{S}(g)(F)(\xi )=e^{-i\frac{u}{4}}e^{\frac{i}{2}%
\sum_{j=1}^{k}x_{j}y_{j}}e^{i\sum_{j=1}^{k}y_{j}\xi _{j}}F(\xi +x),
\label{19}
\end{equation}

\noindent\ and

\[
\Pi _{S}(L)=H=\sum_{j=1}^{k}(-\frac{\partial ^{2}}{\partial \xi _{j}^{2}}%
+\xi _{j}^{2}) 
\]

\noindent is the harmonic oscillator. If $k=1,$ an o.n. basis of
eigenvectors of $H$ in $L^{2}(\mathbb{R})$ is the sequence of Hermite
functions $F_{\mu },\mu \in \mathbb{N}.$ The so called special Hermite
functions \cite[p. 18-19]{T} are, for $\mu ,\mu ^{\prime }\in \mathbb{N}$
and $\varepsilon _{\mu ,\mu ^{\prime }}=sgn(\mu ^{\prime }-\mu ),$

\begin{eqnarray}
\left\langle \Pi _{S}(x,y,0)(F_{\mu }),F_{\mu ^{\prime }}\right\rangle
&=&\Phi _{\mu ,\mu ^{\prime }}(x,y)=\int_{\mathbb{R}}e^{iy\xi }F_{\mu }(\xi +%
\frac{x}{2})F_{\mu ^{\prime }}(\xi -\frac{x}{2})d\xi  \nonumber \\
&=&r_{\mu ,\mu ^{\prime }}(x^{2}+y^{2})e^{-\frac{1}{2}(x^{2}+y^{2})}(x+i%
\varepsilon _{\mu ,\mu ^{\prime }}y)^{\left\vert \mu -\mu ^{\prime
}\right\vert },  \label{32}
\end{eqnarray}

\noindent where $r_{\mu ,\mu ^{\prime }}=r_{\mu ^{\prime },\mu }$ is a one
variable polynomial with real coefficients.

\noindent An o.n basis of eigenvectors of $H$ in $L^{2}(\mathbb{R}^{k})$ is
the sequence $\left( \prod\limits_{j=1}^{k}F_{\mu _{j}}(\xi _{j})\right)
_{\mu \in \mathbb{N}^{k}}$, which gives, for $\mu ,\mu ^{\prime }\in \mathbb{%
N}^{k}$ and $g=(x,y,u),$

\[
\varphi ^{\Pi _{S},\mu ,\mu ^{\prime }}(g)=\left\langle \Pi (g)(F_{\mu
}),F_{\mu ^{\prime }}\right\rangle =e^{-i\frac{u}{4}}\prod\limits_{j=1}^{k}%
\Phi _{\mu _{j},\mu _{j}^{\prime }}(x_{j},y_{j}). 
\]

\noindent By (\ref{32}) the function $z\rightarrow \varphi ^{\Pi _{S},\mu
,\mu ^{\prime }}(zx,zy,z^{2}u)$ is holomorphic on $\mathbb{C}$. Let

\[
R_{a,\delta }=\{\alpha +i\beta \mid \left\vert \alpha \right\vert
<a,\left\vert \beta \right\vert <\delta \}\subset \mathbb{C}. 
\]

\noindent For some constant $C_{a,\delta },$ and $z\in $ $\overline{%
R_{a,\delta }},$

\[
\left\vert \varphi ^{\Pi _{S},\mu ,\mu ^{\prime }}(zx,zy,z^{2}u)\right\vert
\leq C_{a,\delta }e^{\frac{1}{2}a\delta \left\vert u\right\vert
}\prod\limits_{j=1}^{k}e^{\delta ^{2}(x_{j}^{2}+y_{j}^{2})}. 
\]

\noindent We now look for conditions on $a,\delta $ ensuring that the right
hand side lies in $L^{q}(pdg).$ We recall \cite{Hu} that 
\[
p(x,y,u)=\int_{\mathbb{R}}e^{i\lambda u}Q(x,y,\lambda )d\lambda =c_{k}\int_{%
\mathbb{R}}e^{i\lambda u}\prod\limits_{j=1}^{k}\frac{2\lambda }{sh2\lambda }%
e^{-\frac{\lambda }{th2\lambda }(x_{j}^{2}+y_{j}^{2})}d\lambda . 
\]

\noindent Noting that $Q(x,y,\lambda )=\prod\limits_{j=1}^{k}$ $%
Q_{1}(x_{j},y_{j},\lambda )$ is even w.r. to $\lambda ,$ we get, for $q\geq
1,$

\[
\frac{1}{2}\int_{\mathbb{R}}e^{\frac{q}{2}a\delta \left\vert u\right\vert
}p(x,y,u)du\leq \int_{\mathbb{R}}ch(\frac{q}{2}a\delta
u)p(x,y,u)du=Q(x,y,ia\delta \frac{q}{2}). 
\]

\noindent We need the convergence of

\[
\int_{\mathbb{R}^{2k}}\prod\limits_{j=1}^{k}e^{q\delta
^{2}(x_{j}^{2}+y_{j}^{2})}Q_{1}(x_{j},y_{j},i\frac{q}{2}a\delta
)dxdy=c\prod\limits_{j=1}^{k}\int_{\mathbb{R}^{2}}e^{(q\delta ^{2}-\frac{1}{2%
}\frac{qa\delta }{tgqa\delta })(x_{j}^{2}+y_{j}^{2})}dx_{j}dy_{j}, 
\]

\noindent which holds for $qa\delta \leq \frac{\pi }{4}$ and $a>2\delta .$
Thus, taking $a=N\in \mathbb{N},$ $\varphi ^{\Pi _{S},\mu ,\mu ^{\prime }}$
satisfies the assumptions of step 1 on

\[
\Omega =\cup _{N\geq 2}R_{N,\frac{\pi }{4qN}}. 
\]

\noindent Plancherel formula for $\mathbb{H}_{k}$ (see e.g. \cite[Theorem
1.3.1]{T} or \cite[p.154]{CG}) involves the representations

\[
\rho _{h}(x,y,u)=e^{-\frac{i}{4}hu}\Pi _{S}(x,hy,0). 
\]

\noindent By the Stone-Von Neumann theorem \cite[theorem 1.2.1]{T} every
irreducible unitary representation $\Pi $ of $\mathbb{H}_{k}$ satisfying $%
\Pi (0,0,u)=e^{-\frac{i}{4}hu}$ for a real $h\neq 0$ is unitarily equivalent
to $\rho _{h}.$ Hence $\rho _{\beta ^{2}}$ (resp. $\rho _{-\beta ^{2}})$ is
unitarily equivalent to $\Pi _{S}\circ \delta _{\beta }$, (resp. $\Pi
_{S}\circ \sigma \circ \delta _{\beta })$, $\beta >0,$ where $\sigma $ is
the automorphism of $\mathbb{H}_{k}$ defined by $\sigma (x,y,u)=(x,-y,-u)).$

\noindent\ Since $\Pi _{S}(L)=\Pi _{S}\circ \sigma (L),$ we get $\varphi
^{\Pi _{S}\circ \sigma ,\mu ,\mu ^{\prime }}=\varphi ^{\Pi _{S},\mu ,\mu
^{\prime }}\circ \sigma =\overline{\varphi ^{\Pi ,\mu ,\mu ^{\prime }}},$
hence

\[
\mathcal{F}=\{\varphi ^{\Pi _{S},\mu ,\mu ^{\prime }}\circ \delta _{\beta },%
\overline{\varphi ^{\Pi _{S},\mu ,\mu ^{\prime }}}\circ \delta _{\beta
},\beta >0,\mu ,\mu ^{\prime }\in \mathbb{N}^{k}\}. 
\]

\noindent The conditions of step 1 are satisfied by $\varphi ^{\Pi ,\mu ,\mu
^{\prime }}\circ \delta _{\beta },$ replacing $R_{a,\delta }$ by $R_{\beta
a,\beta \delta },$ which ends the proof of theorem \ref{9} for $\mathbb{H}%
_{k}$. Taking remark 5 into account, the set $\cup _{\mu ,\mu ^{\prime
},n}\{h_{n}^{\Pi _{S},\mu ,\mu ^{\prime }},\overline{h_{n}^{\Pi _{S},\mu
,\mu ^{\prime }}}\}$ is total in $L^{2}(\mathbb{H}_{k},pdg).$

\noindent \textbf{Step 3. }We first recall some more facts on
representations and compute the set $\mathcal{F}$ for step 2 stratified
groups. We shall follow Cygan's scheme \cite{Cy}.

\noindent\ \ \ \ \ Let $l\in \mathcal{G}^{\ast }.$ Among the Lie subalgebras 
$\mathcal{M\subset G}$ satisfying $\left\langle l,[X,Y]\right\rangle =0$ for
every $X,Y\in \mathcal{M},$ some have minimal codimension $m_{l}$ and are
denoted by $\mathcal{M}_{l}$. Then the map

\begin{equation}
Z\in \mathcal{M}_{l}\rightarrow e^{i\left\langle l,Z\right\rangle }
\label{88}
\end{equation}

\noindent is a representation of the subgroup $\exp \mathcal{M}_{l}$ and
induces an irreducible unitary representation of $G$ as follows \cite[%
theorems 1.3.3, 2.2.1 and p 41]{CG} : One chooses independent vectors $%
(X_{j})_{i=1}^{m_{l}}$ such that $\mathcal{G}=\mathcal{M}_{l}+span%
\{(X_{j})_{i=1}^{m_{l}}\}.$ For ($g,\xi )\in G$ $\times \mathbb{R}^{m_{l}}$
there exist ($\xi ^{\prime },M)\in \mathbb{R}^{m_{l}}\times \mathcal{M}_{l}$
such that

\[
\exp (\sum\limits_{i=1}^{m_{l}}\xi _{i}X_{i}).\;g=\exp M\;.\exp
(\sum\limits_{i=1}^{m}{}\xi _{i}^{\prime }X_{i}). 
\]

\noindent Then, for $F\in L^{2}(\mathbb{R}^{m_{l}}),$

\begin{equation}
\Pi _{l}(g)(F)(\xi )=e^{i\left\langle l,M\right\rangle }F(\xi ^{\prime }).
\label{29}
\end{equation}

\noindent The set of $\mathcal{C}^{\infty }$ vectors for $\Pi _{l}$ is $%
\mathcal{S(\mathbb{R}}^{k}\mathcal{)}$ \cite[corollary 4.1.2]{CG}. Every
irreducible unitary representation of $G$ is equivalent to a representation
constructed in this way; different $\mathcal{M}_{l},\mathcal{M}_{l}^{\prime
} $ and different $l,l^{\prime }$ \ in the same coadjoint orbit induce
equivalent representations \cite[theorems 2.2.2, 2.2.3, 2.2.4]{CG}.

By Kirillov theory there is an integer $k$ and a set $\mathcal{U}_{0}%
\mathcal{\subset G}^{\ast }$ \ of "generic" orbits with maximal dimension $%
2k,$ such that $m_{l}=k$ for $l\in \mathcal{U}_{0}.$ The Plancherel measure
is supported by $\mathcal{U}_{0}$ \cite[theorem 4.3.10]{CG}.

We now compute such a $\Pi _{l}$ when $G$\textbf{\ }is\textbf{\ }step 2%
\textbf{. }Let $U_{1},..,U_{d}$ be a basis of the central layer $\mathcal{Z}$
and let $\chi _{1},.,\chi _{n}$ be a basis of the first layer $V_{1}$ of $%
\mathcal{G}$.

Let $l\mathcal{\in G}^{\ast }$ and let $\lambda =\sum_{j=1}^{d}\lambda
_{j}U_{j}^{\ast }$ \ be its central part, identified with a vector $\lambda
\in $ $\mathbb{R}^{d}.$ Let $A_{\lambda }$ be the $n\times n$ matrix with
coefficients $\left\langle \lambda ,[\chi _{j},\chi _{h}]\right\rangle .$

By Campbell-Hausdorff formula, for $Y\in \mathcal{G},X\in V_{1},U\in 
\mathcal{Z},g=\exp (X+U),$

\[
\exp Adg(Y)=g\;\exp Y\;g^{-1}=e^{[X,Y]}\exp Y=\exp (Y+[X,Y]), 
\]

\noindent hence the coadjoint orbit of $l,$ \ i.e. $\{l\circ Adg,\;g\in 
\mathcal{G}\}\subset \mathcal{G}^{\ast },$ is $l+range\;A_{\lambda }.$

\noindent We now assume that $l$ lies in $\mathcal{U}_{0},$ so that the
range of $A_{\lambda }$ has dimension $2k.$There exists an orthogonal matrix 
$\Omega _{\lambda }$ such that

\[
A_{\lambda }=\Omega _{\lambda }A_{\lambda }^{\prime }\Omega _{\lambda
}^{\ast } 
\]

\noindent where $A_{\lambda }^{\prime }$ is block diagonal, the non zero
blocks having the form

\begin{equation}
\nu _{j}(\lambda )\left( 
\begin{array}{cc}
0 & 1 \\ 
-1 & 0%
\end{array}%
\right) ,\;\nu _{j}(\lambda )>0.  \label{258}
\end{equation}

\noindent The new basis of $V_{1}$ (defined by the columns of $\Omega
_{\lambda })$ is denoted by $X_{1},Y_{1},..,$ \ \ $%
X_{k},Y_{k},S_{1},...,S_{n-2k},$ so that

\begin{equation}
\left\langle \lambda ,[X_{j},X_{h}]\right\rangle =0=\left\langle \lambda
,[Y_{j},Y_{h}]\right\rangle ,\;\left\langle \lambda
,[X_{j},Y_{h}]\right\rangle =\nu _{j}(\lambda )\delta _{jh},\;1\leq j,h\leq
k.  \label{256}
\end{equation}

\noindent We denote $t=\Omega _{\lambda }(x,y,s)\in \mathbb{R}^{n},$ where 
\[
\sum_{j=1}^{n}t_{j}\chi
_{j}=\sum_{j=1}^{k}x_{j}X_{j}+y_{j}Y_{j}+\sum_{h=1}^{n-2k}s_{h}S_{h}=X+Y+S%
\in V_{1}. 
\]

\noindent Choosing $\mathcal{M}_{l}=\mathcal{Z+}span\{Y_{j},S_{h}\}_{1\leq
j\leq k,1\leq h\leq n-2k}$, let us compute $\Pi _{l}.$ By definition $\Pi
_{l}(\exp u_{j}U_{j})=e^{iu_{j}\lambda _{j}}$. For $g=\exp (X+Z)$ and $%
Z=Y+S, $

\begin{eqnarray*}
\exp (\sum_{j=1}^{k}\xi _{j}X_{j})g &=&\exp [\sum_{j=1}^{k}\xi
_{j}X_{j},X+Z]\;g\;\exp (\sum_{j=1}^{k}\xi _{j}X_{j}) \\
&=&\exp ([\sum_{j=1}^{k}\xi _{j}X_{j},X+Z]+\frac{1}{2}[X,Z])\exp Z\exp X\exp
(\sum_{j=1}^{k}\xi _{j}X_{j}) \\
&=&\exp M\exp (\sum_{j=1}^{k}(\xi _{j}+x_{j})X_{j}).
\end{eqnarray*}

\noindent Hence, by (\ref{29}) and (\ref{256}), for $F\in L^{2}(\mathbb{R}%
^{k}),$

\begin{equation}
\Pi _{l}(g)(F)(\xi )=e^{i\left\langle l,M\right\rangle }F(\xi
+x)=e^{i\sum_{j=1}^{k}\nu _{j}y_{j}(\xi _{j}+\frac{1}{2}x_{j})}e^{i\left%
\langle l,Y+S\right\rangle }F(\xi +x).  \label{31}
\end{equation}

\noindent Since we may replace $l$ by $l^{\prime }$ in the orbit of $l,$ we
may suppose $\left\langle l,Y_{j}\right\rangle =0,$ $1\leq j\leq k.$ In
particular, by (\ref{31}),

\[
\Pi _{l}(X_{j})=\frac{\partial }{\partial \xi _{j}},\Pi _{l}(Y_{j})=i\nu
_{j}\xi _{j},\;1\leq j\leq k,\Pi _{l}(S_{h})=i\left\langle
l,S_{h}\right\rangle I,\;1\leq h\leq n-2k. 
\]

\noindent Since $\Omega _{\lambda }$ is orthogonal, $-L=%
\sum_{j=1}^{k}(X_{j}^{2}+Y_{j}^{2})+\sum_{h=1}^{n-2k}S_{h}^{2},$ which
entails

\[
\Pi _{l}(L)=\sum_{j=1}^{k}-\frac{\partial ^{2}}{\partial \xi _{j}^{2}}+\nu
_{j}^{2}\xi _{j}^{2}+\sum_{h=1}^{n-2k}\left\langle l,S_{h}\right\rangle
^{2}I. 
\]

\noindent A basis of eigenvectors of $\Pi _{l}(L)$ is thus $\left(
\prod\limits_{j=1}^{k}F_{\mu _{j}}(\sqrt{\nu _{j}}\xi _{j})\right) _{\mu \in 
\mathbb{N}^{k}}$. By (\ref{31}) and (\ref{32}), for $g=(x,y,s,u),$

\[
\varphi ^{\Pi _{l},\mu ,\mu ^{\prime }}(g)=e^{i\left\langle \lambda
,u\right\rangle }e^{i\sum_{h=1}^{n-2k}s_{h}\left\langle l,S_{h}\right\rangle
}\prod\limits_{j=1}^{k}\frac{1}{\sqrt{\nu _{j}}}\Phi _{\mu _{j,}\mu
_{j}^{\prime }}(\sqrt{\nu _{j}}x_{j},\sqrt{\nu _{j}}y_{j}). 
\]

\noindent

\noindent Hence, for $z\in R_{a,\delta }$ and some constant $C_{a,\delta },$
with $t=$ $\Omega _{\lambda }(x,y,s),$

\begin{eqnarray*}
\left\vert \varphi ^{\Pi _{l},\mu ,\mu ^{\prime }}(zt,z^{2}u)\right\vert
&\leq &C_{a,\delta }e^{2a\delta \left\vert \left\langle \lambda
,u\right\rangle \right\vert }e^{\delta \sum_{h=1}^{n-2k}\left\vert
s_{h}\left\langle l,S_{h}\right\rangle \right\vert }\prod\limits_{j=1}^{k}%
\frac{1}{\sqrt{\nu _{j}}}e^{\delta ^{2}\nu _{j}(x_{j}^{2}+y_{j}^{2})} \\
&=&e^{2a\delta \left\vert \left\langle \lambda ,u\right\rangle \right\vert
}w_{a,\delta ,l}(x,y,s).
\end{eqnarray*}

\noindent By \cite[corollary 5.5]{Cy} the heat kernel $p(t,u)$ is the
Fourier transform of $CQ(t,\lambda )$ w.r. to the central variables, where

\[
Q(t,\lambda )=\prod\limits_{h=1}^{n-2k}e^{-\frac{1}{2}s_{h}^{2}}\prod%
\limits_{j=1}^{k}Q_{1}(x_{j},y_{j},\frac{\nu _{j}}{4})=Q(t,-\lambda ). 
\]

\noindent Again, we need the convergence \ of

\[
\int_{\mathbb{R}^{n}}w_{a,\delta ,l}^{q}(x,y,s)\prod\limits_{h=1}^{n-2k}e^{-%
\frac{1}{2}s_{h}^{2}}\prod\limits_{j=1}^{k}Q_{1}(x_{j},y_{j},\frac{iqa\delta
\nu _{j}}{2})dxdyds, 
\]

\noindent which holds if $qa\delta \max \nu _{j}\leq \frac{\pi }{4}$ and $%
a>2\delta .$ This ends the proof of theorem \ref{9}.$\blacksquare $

\bigskip \bigskip

Affiliation: \'{e}quipe AGM-UMR 8088

Address: Dept Maths, PST, Universit\'{e} de Cergy, 2 Av. A. Chauvin, 95803,
Cergy, France

E-mail: francoise.piquard@math.u-cergy.fr

\end{document}